\documentclass[twoside,a4paper,11pt]{article}
\newcommand{\heuteIst}{August 10, 2011 }
\usepackage{amsmath,amsfonts,latexsym,amscd,amssymb,amsthm}

\swapnumbers
\theoremstyle{plain}
\newtheorem{theorem}{Theorem}[section]
\newtheorem{lemma}[theorem]{Lemma}
\newtheorem{corollary}[theorem]{Corollary}
\newtheorem{proposition}[theorem]{Proposition}

\theoremstyle{definition}
\newtheorem{definition}[theorem]{Definition}
\newtheorem{example}[theorem]{Example}

\theoremstyle{remark}
\newtheorem{remark}[theorem]{Remark}

\newcommand{\reals}{\mathbb{R}}

\newcommand{\complexs}{\mathbb{C}}

\newcommand{\naturals}{\mathbb{N}}

\newcommand{\integers}{\mathbb{Z}}

\newcommand{\rationals}{\mathbb{Q}}

\DeclareMathOperator{\id}{id}

\newcommand{\abs}[1]{\left\lvert#1\right\rvert} 

\newcommand{\tensor}{\otimes}

\newcommand{\iso}{\cong}
\newcommand{\disjointunion}{\amalg}
\newcommand{\subgroup}{<}
\newcommand{\normalsubgroup}{\lhd}

\DeclareMathOperator{\Aut}{Aut}
\DeclareMathOperator{\im}{im}      

\DeclareMathOperator{\tr}{tr}
\DeclareMathOperator{\pr}{pr}
\DeclareMathOperator{\coker}{coker}

\newcommand{\forget}[1]{}

\newcommand{\innerprod}[1]{\langle #1 \rangle}

{\catcode`@=11\global\let\c@equation=\c@theorem}

\allowdisplaybreaks[2]

\newcommand{\LinnelsC}{\mathcal{C}}
\newcommand{\extendedC}{{{\mathcal{D}}}}

\newcommand{\NeumannN}{\mathcal{N}}
\newcommand{\universalU}{\mathcal{U}}
\newcommand{\RAgroups}{\mathcal{G}}
\newcommand{\amenableGroups}{\mathcal{A}}
\newcommand{\virtabelianGroups}{\mathcal{B}}
\newcommand{\finiteGroups}{(\mathcal{FIN})}
\newcommand{\freegroups}{(\mathcal{FREE})}

\begin{document}
\date{Last compiled \today; last edited  \heuteIst or later}

\title{Integrality of $L^2$-Betti numbers}
\author{
Thomas Schick\thanks{
 e-mail: schick@math-math.gwdg.de
 \protect\\
www:~http://www.uni-math.gwdg.de/schickt/
}\\
Mathematisches Institut --- Georg-August Universit\"at G\"ottingen\\
Germany}

        
\maketitle

\begin{abstract}
The Atiyah conjecture predicts that the $L^2$-Betti numbers of a finite
$CW$-complex with torsion-free fundamental group are integers.
We show that the Atiyah conjecture holds (with an additional
technical condition)
for direct and inverse limits of directed systems of groups for which it is true. As a
corollary it holds for residually torsion-free solvable groups,
e.g.~for pure braid groups or for
positive $1$-relator groups with torsion free abelianization.

Putting everything together we establish a new class of
groups for which the
Atiyah conjecture holds, which contains all free
groups and in particular is closed under taking subgroups, direct
sums,
free products, 
extensions with
elementary amenable quotient,
and under direct and inverse limits of directed systems.

Please take the errata to Schick: ``$L^2$-determinant class and approximation of
$L^2$-Betti numbers'' into account, which are added, rectifying some unproved
statements about ``amenable 
extension''. As a consequence, throughout, amenable extensions should be
extensions with \emph{normal} subgroups.

(non-commutative measure theory)
\end{abstract}

\section{Introduction}\label{sec:introduction}

\begin{remark}\label{rem:gap}
  This is a corrected version of an older paper with the same title.
  The proof of one of the basic results of the earlier version
  contains a gap, as was kindly pointed out to me by Pere Ara.
  This gap could not be fixed. Consequently, in this new version
  everything based on this result had to be removed.
\end{remark}

In \cite{Atiyah(1976)} Atiyah introduced $L^2$-Betti numbers of closed
manifolds in terms of the kernel of the Laplacian on the universal
covering. There, he asked what the possible values of these numbers
are, in particular whether they are always integers if the fundamental
group is torsion-free.  We call this the Atiyah conjecture. More
precisely, we consider the following algebraic question:

\begin{definition}\label{AtiyahCon} Let $\integers\subset
  \Lambda\subset \rationals$ be an additive subgroup
  of the rationals, and $K$ a subring of $\complexs$.
  We say a discrete group $G$ fulfills the
 \emph{Atiyah
conjecture of order $\Lambda$ over $KG$} if 
\begin{equation*}
              \dim_G (\ker A)\in \Lambda \qquad\forall A\in
M(m\times n, K G),
              \end{equation*}
where $\ker A$ is the kernel of the induced map $A\colon l^2(G)^n\to
l^2(G)^m$. Let $\pr_{\ker A}$ be the projection onto $\ker A$. Then
\begin{equation*}
\dim_G(\ker A):=\tr_G(\pr_{\ker A}):=\sum_{i=1}^n \innerprod{\pr_{\ker A}
  e_i,e_i}_{l^2(G)^n},
\end{equation*}
where $e_i\in l^2(G)^n$ is the vector with the trivial element of
$G\subset l^2(G)$
at the $i^{\text{th}}$-position and zeros elsewhere. $\tr_G$ is the
canonical finite
trace on the von Neumann algebra of operators on $l^2(G)^n$ commuting
with the right $G$-action, the so called \emph{Hilbert $\NeumannN
  G$-module maps}.

The group $G$ is said to fulfill the \emph{strong Atiyah conjecture over
 $KG$} if $\Lambda$ is
 generated by $\{\abs{F}^{-1}|\; F\subgroup G \text{ finite}\}$.
\end{definition}

\begin{remark}
  In \cite{Grigorchuk-Linnell-Schick-Zuk(2000)} it is shown that the so
  called Lamplighter group $G$ does not satisfy the Atiyah conjecture over
  $\integers
  G$. However, there is no bound on the orders of finite subgroups of
  $G$, in contrast to all examples of groups for which the strong
  Atiyah conjecture is known.
\end{remark}

In Definition \ref{AtiyahCon}, we can replace the kernels  by the closures of the
images or by the cokernels because these numbers by additivity of the
$G$-dimension \cite[1.4]{Lueck(1997)} pairwise sum up to
the dimension of the
domain or range, which is an integer. Replacing $A$ by $A^*A$ we may also
assume that $n=m$. Because we can multiply $A$ with a unit in
$Q$ without changing its kernel, the Atiyah conjecture over
$KG$ implies the Atiyah conjecture over $QG$, where $Q$ is the quotient field of
$K$ in $\complexs$. Usually therefore we will assume that
$K$ is a field.

The Atiyah conjecture of order $\Lambda$ over $\integers G$ is
equivalent to the statement that
all $L^2$-Betti
numbers of finite
$CW$-complexes with fundamental group $G$ are elements of $\Lambda$
\cite[Lemma 2.2]{Lueck(1997)}, i.e.~to the original question of Atiyah.

Fix a coefficient field $K\subset \complexs$. If not stated
otherwise, the term ``Atiyah conjecture'' means ``Atiyah conjecture
over $KG$''.



We introduce a new method of proof for the Atiyah conjecture
which is based on the approximation results in
\cite{Schick(1998a)}. In \cite{Schick(1998a)} it is proved that in
many cases if a group $G$ is the
direct or inverse limit of a directed system of groups $G_i$, then
$L^2$-Betti numbers over $G$ are limits of $L^2$-Betti numbers defined
over the groups $G_i$. More precisely:
\begin{definition}\label{defG}  (compare
  \cite[1.11]{Schick(1998a)})
  Let $\RAgroups$ be the smallest class of groups which contains the
  trivial group, is subgroup closed and is closed under 
  direct or inverse
  limits of directed systems, and under extensions with amenable quotient.

The amenable quotient is not
necessarily a quotient 
group, one  relaxes this to a suitable notion of an amenable
(discrete) quotient space (by a subgroup which is not normal) (compare
\cite[4.1]{Schick(1998a)}).
\end{definition}
This class of groups is very large. It contains most groups which
naturally occur in geometry, e.g.~all residually finite groups. It
would be interesting to find an example of a group which does not
belong to $\RAgroups$.

\begin{proposition}\label{tflimits}
  Let $G_i\in\RAgroups$ (with $\RAgroups$ as in Definition \ref{defG})
  be a directed system of torsion-free groups
  which fulfill the strong Atiyah conjecture. Then their (direct or
  inverse) limit $G$ fulfills the strong Atiyah conjecture over $\rationals G$.
\end{proposition}

\begin{corollary}
  The strong Atiyah conjecture over $\rationals G$ is true if $G$ is a 
  pure braid group, or a
  positive $1$-relator group (in the sense of \cite{Baumslag(1971)})
  with torsion free abelianization.
\end{corollary}
The corollary will be explained in Example \ref{newex}.

 There
is one significant difference to all other
results about the Atiyah conjecture obtained so far: the approximation
result is proved only over $\rationals G$ instead of $\complexs
G$. This does not effect the original question of Atiyah. But it is
relevant for the zero divisor conjecture (by Lemma \ref{zerodiv}). In a
forthcoming paper \cite{Linnell-Schick(2000)} we will show how to
enlarge $\rationals$ at least to the field of algebraic numbers in $\complexs$.

Proposition \ref{tflimits} can be used to give an alternative proof of
the Atiyah
conjecture for free groups, compare Example \ref{freeex}.

Linnell obtained the  most general positive results about the Atiyah
conjecture so far \cite{Linnel(1993),Linnell(1998)}.
His interest stems from
the zero divisor conjecture and the following observation
(compare e.g.~\cite{Linnel(1992)}):
\begin{lemma}\label{zerodiv}
  If for a torsion-free group the strong Atiyah conjecture is true
  over $K G$ then the ring $K G$ has no non-trivial zero divisors.
\end{lemma}

 Linnell uses essentially two lines of argument for the Atiyah conjecture:
algebraic considerations using $K$-theoretic information (in
particular Moody's induction theorem \cite[Theorem 1]{Moody(1989)})
to obtain the Atiyah conjecture for elementary amenable groups or for
extensions with elementary amenable quotient.

For a free group $F$, Linnell devised a completely different argument,
making use of Fredholm module techniques which were used before
to prove that there are no non-trivial projectors in $C^*_rF$. 

Linnell defines \cite{Linnel(1993),Linnell(1998)}:
\begin{definition}
  \label{firstdefC}
 The class $\LinnelsC$ is the smallest class of groups which contains all
  free groups and is closed under directed union and
  extensions with virtually abelian  quotients. 
The class $\LinnelsC'$ consists of extensions of
  direct sums of free groups with elementary amenable quotient.
\end{definition}

Putting his two approaches together, Linnell proves in
\cite[1.5]{Linnel(1993)}
\begin{theorem}\label{LinnellsTheorem}
  \label{LCprime} The strong Atiyah conjecture over $\complexs G$ is
  true for groups  in the
classes $\LinnelsC$ 
which have a bound on the orders of finite subgroups.
\end{theorem}

The proof of the Atiyah conjecture for the class $\LinnelsC'$ given in
\cite{Linnell(1998)} contains a gap (which parallels the gap mentioned 
in Remark \ref{rem:gap}). However, we will shown in Corollary
\ref{corol:LinnelsCprime} how to fill this gap and conclude that the
strong Atiyah conjecture over
$\rationals G$ is true if $G\in\LinnelsC'$.


We use the new results to define the following
class of groups for which the Atiyah conjecture is true:
\begin{definition}\label{defextC}
  Let $\extendedC$ be the smallest non-empty class of groups such that:
  \begin{enumerate}
  \item \label{aex} If $G$ is torsion-free and $A$ is elementary
    amenable, and we have a projection
    $p\colon G\to A$ such that $p^{-1}(E)\in\extendedC$ for every
    finite subgroup $E$ of $A$, then $G\in\extendedC$.
   \item \label{sgr} $\extendedC$ is subgroup closed.
   \item\label{lim} Let $G_i\in \extendedC$ be a
  directed system of groups
    and $G$ its (direct or inverse) limit. Then $G\in\extendedC$.
\end{enumerate}
\end{definition}

Observe that $\extendedC$ contains
only torsion-free groups.

\begin{theorem}\label{strongAt}
  If  $G\in\extendedC$ then the strong Atiyah
  conjecture over 
  $\rationals G$ is true for $G$.
\end{theorem}

\begin{corollary}\label{corol:zero_divisors_in_D}
  If $G\in\extendedC$ then the ring $\rationals G$ does not contain
  non-trivial zero divisors.
\end{corollary}

\begin{proposition}\label{classprop}
  All torsion-free groups in $\LinnelsC$ and $\LinnelsC'$ are also
  contained in $\extendedC$. Moreover,
 $\extendedC$ is
  closed under direct sums, direct products and free
  products.
%
\end{proposition}

$\LinnelsC$  is not
closed under direct sums. $\LinnelsC'$ has this property, but it is
not closed under free product. $\extendedC$ is an enlargement which
repairs these deficits.

\begin{example}
  It follows that
  $G:=\integers*(\integers*\integers \times
  \integers*\integers)$ fulfills the Atiyah conjecture over
  $\rationals G$. The $L^2$-Betti numbers of $G$ are computed as
  $b_1^{(2)}(G) =1 =
  b_2^{(2)}(G)$ and $b_k^{(2)}(G)=0$ for $k\ne 1,2$. By 
  \cite[9.3]{Reich(1999)} $G\notin
  \LinnelsC$. Using ideas of Reich, one can also show $G\notin
  \LinnelsC'$,
 so that indeed we cover additional groups.
\end{example}

\subsection*{Organization of the paper}
Section
\ref{sec:limits} gives an overview over elementary results about the
Atiyah conjecture. We also prove Proposition \ref{tflimits} and
discuss applications of this proposition. In Sections
\ref{secextension} and \ref{sec:proofs} we prove Theorem
\ref{strongAt} and Proposition \ref{classprop}. Section \ref{sec:gen}
describes possible 
generalizations. In Section
\ref{treechar} we extend Linnell's proof of the Atiyah conjecture for
free groups to another class of groups, the so called treelike groups, 
and characterize the treelike groups.

\section{The Atiyah conjecture for limits of groups}\label{sec:limits}

Here we collect a few statements where the
 Atiyah conjecture for some groups implies its validity for other
 groups. 

The following propositions are well known.
\begin{proposition}\label{finext} (\cite[8.6]{Linnell(1998)})\\
  If $H$ is a subgroup of index $n$ in $G$ and $H$ fulfills the
Atiyah conjecture of order $\Lambda\subset\rationals$ over $KH$, then
$G$ fulfills the Atiyah
conjecture of order $\frac{1}{n}\Lambda$ over $KG$.
\end{proposition}

\begin{proposition}\label{subgroups}
  If $G$ fulfills the  Atiyah conjecture of order
  $\Lambda\subset\rationals$ over $KG$, and if $U$ is a
  subgroup of $G$, then
$U$ fulfills the Atiyah conjecture of order $\Lambda$ over $KH$.
\end{proposition}
\begin{proof}
  This follows from the fact that the $U$-dimension of the kernel of
  a matrix over $K U$ acting on $l^2(U)^n$ coincides with the
  $G$-dimension of the
  same matrix, considered as an operator on $l^2(G)^n$
  \cite[3.1]{Schick(1998a)}.
\end{proof}

\begin{proposition}  \label{unions}
  Let $G$ be the directed union of groups $\{G_i\}_{i\in I}$ and
  assume that each
  $G_i$ fulfills the Atiyah conjecture of order $\Lambda_i\subset
  \rationals$ over $KG_i$. Then $G$ fulfills the Atiyah conjecture of order
  $\Lambda$ over $KG$, where $\Lambda$ is the additive subgroup of
  $\rationals$ generated by $\{\Lambda_i\}_{i\in I}$.
\end{proposition}
\begin{proof}
  A matrix over $KG$, having only finitely many non-trivial
  coefficients, already is a matrix over $KG_i$ for some $i$. The
  $G_i$-dimension and the $G$-dimension of the kernel of the
  matrix coincide by \cite[3.1]{Schick(1998a)}.
\end{proof}

For groups in $\RAgroups$ we have good
approximation results
for $L^2$-Betti numbers \cite[6.9]{Schick(1998a)}. We use
  these to
give a proof of the Atiyah conjecture for suitable groups.
For this it is
necessary to work over the rational
group ring instead of the complex group ring.

\begin{proposition}\label{limits}
  Let $G_i$ be a directed system of groups which fulfill the 
Atiyah conjecture over $\rationals G_i$ of common order
$\frac{1}{L}\integers$. Suppose all the
groups $G_i$ lie in the class $\RAgroups$ (defined in
Definition \ref{defG}). Then their direct or inverse
limit $G$ fulfills the  Atiyah conjecture of order
$\frac{1}{L}\integers$ over $\rationals G$.
\end{proposition}
This is a direct consequence of  the assumption and the approximation result
\cite[6.9]{Schick(1998a)}, since
 there it is shown that
each $L^2$-Betti number over the limit group is the limit of
$L^2$-Betti numbers over the groups in the sequence. By assumption,
all of these lie in $\frac{1}{L}\integers$, which is a closed subset of
$\reals$. Therefore the same is true for the limit $L^2$-Betti number.

Note that this implies the strong Atiyah conjecture if we know that
$L$ is the least
common multiple of the orders of finite subgroups of $G$ (provided
this number exists), in particular if all $G_i$ are torsion-free and
fulfill the 
strong Atiyah conjecture. Hence Proposition \ref{tflimits} is a
special case of Proposition \ref{limits}.

\begin{example}\label{freeex}
  We use Proposition \ref{limits} to give a different proof that a free
  group $F$ fulfills the
  strong Atiyah 
  conjecture over $\rationals F$. By \cite{Magnus(1935)} free groups
  are residually 
  torsion-free nilpotent. More precisely, if we consider the
  descending central series
  $F\supset F_1\supset F_2\supset\dots$ then $\bigcap F_k=\{1\}$ and
  $F/F_k$ is torsion-free. It follows that the groups $F/F_k$ are
  nilpotent and poly-$\integers$. Using the methods of the
  ring-theoretic parts of \cite{Linnel(1993)} it is easy to prove the
  strong Atiyah conjecture (over the complex group ring) for
  poly-$\integers$ groups. We will repeat the proof in Section
  \ref{secextension}. Alternatively, one can rely on a different type
  of ring theory and in particular avoid the use of the rings $DG$ and
  $\universalU G$ of Definition \ref{defunivU}:
 If $G$ is poly-$\integers$, then $\complexs G$ is
  Noetherian \cite[8.2.2]{Rowen(1988b)}. Using Moody's induction theorem
  \cite[Theorem 1]{Moody(1989)} we see that the $G$-theory of
  $\complexs G$ is trivial,
  which by
  \cite[4.6]{Swan(1968)} coincides with the
  $K$-theory. Therefore a finite projective resolution of a $\complexs
  G$-module can be replaced by a finite free resolution. By
  \cite[8.2.18]{Rowen(1988b)} $\complexs G$ has
  finite cohomological dimension. Together this implies that every
  finitely generated
  module over $\complexs G$ has 
  a finite free resolution.
  If $A$ is a matrix
  over $\complexs G$, its
  kernel in $(\complexs G)^n$ is finitely generated and we define
  its $G$-dimension as the alternating sum of the ranks of the modules
  in a finite free resolution. One can prove that this number,
  which is an integer, coincides with the $G$-dimension of $\ker
  A\subset l^2(G)^n$ (this is e.g.~done in the thesis
  \cite[5.4.1]{Bratzler(1997)}).

  We therefore arrive at the inverse system $\dots\to F/F_2\to F/F_1$
  where each of the quotients  $F/F_k$ is torsion-free and fulfills the
  strong Atiyah conjecture and lies in $\RAgroups$. $F$ is a
  subgroup of the inverse
  limit of this system. Application of Proposition \ref{tflimits} and
  Lemma \ref{subgroups} gives the result.

There is of course one drawback:
   the 
  approximation method proves 
  the Atiyah conjecture only for matrices over $\rationals F$ instead of
  $\complexs F$.
\end{example}

In view of Proposition \ref{limits} the following lemma is of
importance.
\begin{lemma}
  \label{LinnelsC_in_RA} 
  The class $\extendedC$ is contained in $\RAgroups$.
\end{lemma}
\begin{proof}

By definition, $\RAgroups$ is subgroup
  closed and closed under
  direct and inverse limits of directed systems, and under extensions
  with elementary amenable quotients. It follows
  that $\extendedC$ is contained in $\RAgroups$.

\end{proof}

\begin{corollary}\label{resC}
  The strong Atiyah
  conjecture over $\rationals G$ is true if $G$ is residually torsion
  free of class $\LinnelsC$, in particular if $G$ is residually
  torsion-free solvable (or more generally residually torsion-free elementary amenable).
\end{corollary}
\begin{proof}
  Let $\mathcal{X}$ be a class of groups.
    We call a group $G$ residually of class $\mathcal{X}$ if 
  $G$ has a sequence of normal subgroups $G\supset H_1\supset H_2\supset\dots$
  with $\bigcap_{k\in\naturals} H_k =\{1\}$ and such that the
  quotients $G/H_k\in\mathcal{X}$ $\forall k\in\naturals$.
Another (more common)  definition of residuality is to require for every $g\in
  G$ a homomorphism $\phi_g:G\to X_g$ with $\phi_g(g)\ne 1$ where
  $X_g\in\mathcal{X}$. If $\mathcal{X}$ is closed under direct
  sums  then for countable groups $G$ the two
  definitions are equivalent. This is the case for the class of
  torsion-free elementary amenable groups (but not for
  the class $\LinnelsC$).

  If $G$ has a sequence of normal subgroups $H_1\supset
  H_2\supset\dots$ with trivial intersection such that
  $G/H_k\in\LinnelsC$ and $G/H_k$ are torsion-free, then $G$ is a
  subgroup of the inverse limit of
  $(G/H_k)_{k\in\naturals}$. Combining Lemma \ref{LinnelsC_in_RA} and
  Proposition \ref{classprop}, $\LinnelsC\subset\RAgroups$. Because
  of Theorem \ref{LinnellsTheorem}
  we can apply Proposition \ref{tflimits}
  and Lemma \ref{subgroups} to conclude that $G$ fulfills the strong
  Atiyah conjecture over $\rationals G$.
\end{proof}

\begin{example}\label{newex}
  Let $F$
  be a free group on the free generators $\{a_i\}_{i\in \naturals}$
  and assume $c\in F$ generates its
  own centralizer. Choose $V\subset\naturals$. Let $A$ be a
  free abelian countably generated group
  $A$ and $r$  one of the  free generators.
  
The following groups are residually torsion free of class
  $\LinnelsC$ and hence by Corollary \ref{resC} fulfill the strong
  Atiyah conjecture over the rational group ring:
  \begin{enumerate}
\item Pure (also called colored) braid groups.
  \item Positive  $1$-relator groups $G$ with $H_1(G,\integers)$
    torsion free. A group is a $1$-relator group if it has a
    presentation 
    \begin{equation*}
G=<g_1,\dots,g_n| R:=\prod_{k=1}^N g_{i_k}^{n_k}=1>
\end{equation*}
    with one relation $R$. It is called positive if one finds such a
    presentation with $n_k\ge 0$ for all $k$. By \cite[IV.5.2.]{Lyndon-Schupp(1977)}
    $G$ is torsion free if and only
    if $R$ is not a proper power in the free group generated by
    $g_1,\dots,g_n$. It follows from the presentation that
    $H_1(G)=G/[G,G]$ is torsion free if and only if the greatest
    common divisor of $s_1,\dots,s_n$ is $1$, where
    $s_r=\sum_{i_k=r}n_k$ is the exponent sum of $g_r$ in $R$. Observe
    that therefore $H_1(G,\integers)$ torsion free for a positive
    $1$-relator group implies that $G$ is torsion free.
  \item  $
    E_1:=<r,a_i; i\in V\;|\; c=r^n> $;
\item $
    E_2:=<a_i; i=1,\dots,k\;|\; a_{1}^{n_1}a_{2}^{n_2}\dots a_k^{n_k}=1>$, with $k>3$
    and $n_i\in\integers$;
\item $E_3:= F*_{c=r}A$.
\end{enumerate}
\end{example}
\begin{proof}
 First, we deal with the positive $1$-relator groups.
  If $R=1$, $G$ is free and we have nothing to prove.
  By \cite[Theorem 1]{Baumslag(1971)} positive $1$-relator groups are
  residually
  solvable. We want to prove that all the derived series quotients
  $G/G^{(n)}$ are torsion free. We can assume that $R\ne 1$, and since all
  $n_k$ are non-negative, the second differential in the presentation
  $2$-complex $X^{(2)}$ of $G$ is non-trivial (i.e.~at least one of
  the $s_k$ is not zero), and therefore $H_2(X^{(2)},\integers)=0$. By
  assumption $H_1(G,\integers)$ is torsion free. Now
  \cite[Theorem A]{Strebel(1974)} implies that $G/G^{(n)}$ is torsion free for
  all $n$.

  By \cite[Theorem
  2.6]{Falk-Randell(1988)}, the pure braid groups are residually torsion-free
  nilpotent.
  The group $E_1$ is residually torsion-free nilpotent by \cite[Theorem
  1]{Baumslag(1968)}. $E_2$ is residually free by
  \cite[p.~414]{Baumslag(1967)}. $E_3$    is residually free \cite[Theorem
  8]{Baumslag(1967)}. Free groups and therefore $E_1,E_2,E_3$ are
  residually
  torsion-free solvable. Every nilpotent and every solvable group is contained
  in $\LinnelsC$.
\end{proof}

\begin{remark}

  In \cite{Linnell-Schick(2000a)} a condition is given when the Atiyah
  conjecture for a group $H$ implies the Atiyah conjecture for finite
  extensions of $H$. The pure braid groups satisfy this condition, so
  that in particular the strong Atiyah conjecture for the full braid
  groups follows.\\
  The condition is also satisfied for many positive one-relator groups.
\end{remark}


\section{Extensions with elementary amenable quotient}\label{secextension}

The following proposition is implicit in Linnell's work. He  only
deals with the complex group ring, therefore we give a complete
proof here.

\begin{proposition}\label{amext}
  Let $1\to H\to G\to A\to 1$ be an exact sequence of
  groups. Assume that $G$ is torsion free and $A$
  is elementary amenable. For every finite subgroup
  $E\subgroup A$ let $H_E$ be the inverse
  image of $E$ in $G$. Assume that $K=\bar K\subset\complexs$ and for all
  finite subgroups
  $E\subgroup G$ that $H_E$ fulfills
  the strong Atiyah
  conjecture over $KH_E$. Then $G$
  fulfills
  the strong Atiyah conjecture over $KG$.

\end{proposition}
\begin{corollary}\label{tfamext}
 Suppose
  $H$ is torsion-free and fulfills the strong Atiyah conjecture over
  $KH$ with
  $K=\bar K$. If $G$ is an
  extension of $H$ with elementary amenable torsion-free quotient then
  $G$ fulfills the strong Atiyah conjecture.
\end{corollary}
\begin{proof}
  By assumption, the only finite subgroup of $G/H$ is the trivial
  group and the Atiyah conjecture is true for its inverse image $H$.
\end{proof}

We essentially follow Linnell's lines in
  \cite{Linnel(1993)} to prove  Proposition \ref{amext}. First we
  repeat a few lemmas.

\begin{definition}\label{defunivU}
  Given the group von Neumann algebra $\NeumannN G$ we define
  $\universalU G$ to be the algebra of all unbounded operators
  affiliated to $\NeumannN G$ (compare e.g.~\cite[Section 8]{Linnell(1998)}),
  i.e.~such that all their spectral projections belong to $\NeumannN
  G$.

  We have to consider $KG\subset
  \NeumannN G\subset \universalU G$. Let $DG$ be the division closure
  of $KG$ in $\universalU G$, i.e.~the smallest subring of
  $\universalU G$ which contains $KG$ and which has the property that
  whenever $x\in DG$ is invertible in $\universalU G$ then $x^{-1}\in
  DG$.
\end{definition}

 The Atiyah conjecture is strongly connected to ring theoretic
  properties of $DG$:
\begin{lemma}\label{skew}
  Let $G$ be a torsion-free group. $G$ fulfills the strong
  Atiyah conjecture over $KG$
  if and only if the division closure $DG$ of $KG$ in
  $\universalU G$ is a skew field.
\end{lemma}
\begin{proof}
  The ``only if'' part is proved in \cite{Linnell(1998)}. Linnell states the
  result only for $K=\complexs$ but the proof carries over verbatim to
  the more general case. For the converse suppose $DG$ is a skew
  field. We have to show that for $A\in M(m\times n, KG)$ the
  $G$-dimension of
  $\coker(A\tensor_{KG} \id_{l^2G})$ is an integer, or
  equivalently that the $G$-dimension of
  $\coker(A\tensor_{KG}\id_{\universalU G})$ is an
  integer. Since tensor products are right exact, this amounts to the
  fact that the $G$-dimension of
  $(\coker(A)\tensor_{KG}DG)\tensor_{DG}\universalU G$
  is integral. However, $\coker(A)\tensor_{KG}DG$ is a
  (finitely generated, since $A$ is a finite matrix) module over the
  skew field $DG$, therefore isomorphic to $(DG)^N$ for a
  suitable integer $N$. Then
  $\coker(A)\tensor_{KG}\universalU G\cong
  (\universalU G)^N$ which has $G$-dimension $N\in\naturals$.
\end{proof}

 \begin{lemma}\label{polinv}
   Suppose $H$ is a normal subgroup of $G$ with infinite cyclic
   quotient. Let $(\universalU H)G$ be the subring of $\universalU G$
   generated by $G$ and $\universalU H$. If $t\in G-H$ and
   $x=1+q_1t+\dots+q_kt^k\in(\universalU H)G$ with $q_i\in\universalU
   H$ then $x$ is invertible in $\universalU G$ and in particular is
   not a zero divisor.
 \end{lemma}
 \begin{proof}
   If $q\in \universalU H$ then one finds $\alpha,\beta\in\NeumannN H$
   with $q=\alpha\beta^{-1}$ by \cite[Theorem 1 and proof of Theorem
   10]{Berberian(1982)} or \cite[2.7]{Reich(1999)}. Applying the same
   to $q^*$ we may as well achieve $q=\beta^{-1}\alpha$. Apply this now
   inductively to $q_1=\beta_1^{-1}\alpha_1$,
   $\beta_1q_2=\beta_2^{-1}\alpha_2,\dots$ to get the non-zero divisor
   $\beta=\beta_1\dots\beta_k\in\NeumannN H$ (which is invertible in
   $\universalU H$) such that $\beta q_i\in\NeumannN H$ for
   $i=1,\dots,k$. It follows from \cite[Theorem 4]{Linnel(1992)} that
   $\beta x$ is not a zero divisor in $\NeumannN G$, therefore it
   becomes invertible in $\universalU G$. 
 \end{proof}

\begin{lemma}\label{star}
  Let $DG$ be the division closure of $KG$ in $\universalU G$, where
   $K=\bar K$. If $\alpha\in\Aut(G)$ then
  \begin{equation*}
    \alpha(DG)=DG^*=DG .
  \end{equation*}
\end{lemma}
\begin{proof}
  Of course $\alpha(DG)$ is the division closure of $\alpha(KG)$ in
  $\alpha(\universalU G)$. Since $\alpha(KG)=KG$ and
  $\alpha(\universalU G)=\universalU G$ we conclude
  $\alpha(DG)=DG$. The proof for the anti-automorphism $*$ is
  identical. Here $KG=KG^*$ since $K$ is closed under complex conjugation.
\end{proof}

\begin{definition}\label{Def:Groups}
Mainly for proofs by induction we will use the following constructions
of classes of groups in
addition to Definition \ref{defextC} and Definition \ref{defG}:
\begin{itemize}
\item For a class $\mathcal{X}$ of groups, $L\mathcal{X}$ denotes the
  class of all groups which are locally of class $\mathcal{X}$,
  i.e.~every finitely generated subgroup belongs to
  $\mathcal{X}$. Such groups are directed unions of groups
  in $\mathcal{X}$.
\item For two classes of groups $\mathcal{X}$ and $\mathcal{Y}$, the class
  $\mathcal{X}\mathcal{Y}$ consists of all groups $G$ with normal
  subgroup $H\in\mathcal{X}$ and quotient $G/H\in\mathcal{Y}$.
Let $\mathcal{X}^\times$
  consist of all finite direct sums of groups in $\mathcal{X}$.
\item 
$\freegroups$ is the class of free groups, 
 $\finiteGroups$ is
  the class of finite groups, and $\virtabelianGroups$ is the class of finitely generated
   virtually abelian groups.
\item The elementary amenable groups are denoted by $\amenableGroups$.
  This is the smallest class of groups which contains
  all abelian and all finite groups and is closed under extensions and
  directed unions.
\end{itemize}
\end{definition}

\begin{proof}[Proof of Proposition \ref{amext}]\strut\\
\noindent  Case 1: $G/H$  finitely generated free
  abelian.\\
  By induction we can immediately reduced to the case where $G/H$ is
  infinite cyclic.
   Since $DH$ is a skew field an
   arbitrary non-zero element of $DH*G/H$ has after multiplication with
   units in $DH*G/H$ the shape $1+q_1t+q_2t^2+\dots q_kt^k$ where $t$
   is a generator of $G/H$. Lemma \ref{polinv} implies
   therefore that each element of $(DH)G-\{0\}$ becomes invertible in
   $\universalU G$, in particular $(DH)G$ has no zero
   divisors. Therefore its Ore completion is a skew field which embeds
   into $\universalU G$. It embeds into $DG$ and is division closed,
   therefore both coincide, $DG$ is a skew field and by Lemma
   \ref{skew} the strong Atiyah conjecture holds for $G$.

   If, as a next step, $G/H$ is finitely generated free abelian,
   induction immediately implies that the strong Atiyah conjecture
   holds for $G$.

\noindent
  Case 2: $G/H$ finitely generated virtually abelian.\\ It has a finitely
  generated free abelian
  normal subgroup $A_0$. Its inverse image $G_0$ in $G$ fulfills the
  strong Atiyah conjecture by the previous step.
  Since $H$ fulfills the strong Atiyah conjecture,
  by Lemma \ref{skew} the division closure $DH$ of $KH$ in
  $\universalU H$ is a skew
  field. Conjugation with an element of $G$ gives an automorphism
  of $H$. Using Lemma
  \ref{star}, by
  \cite[2.1]{Linnel(1993)} the subring $(DH)G$ generated by $DH$ and
  $G$ is a crossed product: $(DH)G=DH*(G/H)$.  The same applies to every
   finite extension $H\normalsubgroup H_E$, therefore $(DH)H_E =
   DH*(H_E/H)$. By    \cite[4.2]{Linnel(1993)} for these finite extensions
 $(DH)H_E=DH_E$, since $DH$ is a
   skew field and hence Artinian.

   By assumption all the groups $H_E$ are torsion-free and fulfill
  the strong Atiyah conjecture. By Lemma \ref{skew}
  $DH_E=DH*E$ then is a skew field for every finite $E\subgroup G/H$. 
  Because of \cite[4.5]{Linnell(1998)},
   $DH*G/H=(DH)G$ is an Ore domain. The proof shows that its Ore
  completion (a skew field!) is obtained by
  inverting $DH*(G_0/H)-\{0\}$.
   Since $DH*(G_0/H)\subset DG_0$, and the latter is a skew field contained
  in $\universalU G$, all inverses
  $\left(DH*(G_0/H)-\{0\}\right)^{-1}$ are contained in $\universalU
  G$. Therefore the same is true for the Ore completion of $DH*G/H$,
  which is a skew field and as a result coincides with the division
  closure $DG$. By Lemma
  \ref{skew} the strong Atiyah conjecture holds for $G$.

\noindent
Case 3: $A$ elementary amenable.\\
 This is done by
  transfinite induction, where we use roughly the
  description of elementary amenable groups given in \cite[Section
  3]{Kropholler-Linnel-Moody(1988)}: Let  $\amenableGroups_0$ be the class
  of finite groups, and for an ordinal $\alpha$ set
  $\amenableGroups_{\alpha+1}:=(L\amenableGroups_\alpha)\virtabelianGroups$.
  If $\alpha$ is a limit ordinal set 
  $\amenableGroups_{\alpha}=\bigcup_{\beta<\alpha}\amenableGroups_{\beta}$.
  Then $\amenableGroups=\bigcup_{\alpha\ge
 0}\amenableGroups_{\alpha}$. One can easily show, as in
  \cite[4.9]{Linnel(1993)}, that each $\amenableGroups_\alpha$ and
  $L\amenableGroups_\alpha$ is
  subgroup closed and closed under extensions with finite
  quotient.

  The case $A\in\amenableGroups_0$ is trivial, because then the
  conclusion is part of the assumptions.

  Let $\alpha$ be the least ordinal with
  $A\in\amenableGroups_\alpha$. Then $\alpha=\beta+1$. Assume by
  induction that the statement is true for $\amenableGroups_\beta$.
  If $A=\bigcup_{i\in I}A_i$ with $A_i\in\amenableGroups_\beta$, then
  $G=\bigcup G_i$ where $G_i$ is the inverse image of $A_i$ in $G$,
  i.e.~an extension of $H$ by $G_i$. By induction and Lemma
  \ref{unions} the Atiyah conjecture is true for $G$, therefore it
  is true if the quotients belong to $L\amenableGroups_\beta$.

  Let now $A$ be an extension with kernel $B\in L\amenableGroups_{\beta}$ and quotient $Q$ which is
  finitely generated virtually abelian.  We want to apply the
  induction hypothesis and Case 2 to $1\to G_1\to G\to Q\to
  1$, where for a finite subgroup $E$ of $Q$ we let $G_E$ be the
  inverse image of $E$ in $G$ (in particular
  $G_1$ is the inverse image of $B$). To apply the results of Case 2
  we have to establish the Atiyah conjecture for all the groups
  $G_E$. We obtain exact
  sequences $1\to H\to G_E\to B_E\to 1$  where $B_E$ is a finite extension of
  $B$. Hence $B_E\in L\amenableGroups_\beta$. Note that $B_E$ is
  a subgroup of $A$, 
  therefore each finite subgroup of $B_E$ is a finite subgroup of $A$,
  too. This implies that the assumptions of the proposition are
  fulfilled for $1\to H\to G_E\to B_E\to 1$ and therefore $G_E$
  fulfills the strong Atiyah conjecture by induction.

  Hence the
  assumptions of Case 2 are fulfilled for $1\to
  G_1\to
  G\to Q\to 1$ and the result follows.
\end{proof}

In the course of the proof we also established:
\begin{lemma}\label{oreloc}
  If, in the situation of Proposition \ref{amext}, $G/H$ is finitely
  generated virtually
  abelian then the ring $DG$ is an Ore
  localization of $DH*G/H$.
\end{lemma}

\begin{example}
  Let $G\in\RAgroups$ (as defined in  Definition \ref{defG}) be a torsion-free
  group which fulfills the
  strong Atiyah 
conjecture with $K=\bar K\subset\complexs$, and $f:G\to G$ a group
  homomorphism. Then
the mapping torus group 
\begin{equation*}
T_f:=<t,G|\;tgt^{-1}=f(g), gh=(gh) \quad g,h\in G>
\end{equation*}
 fulfills the strong Atiyah conjecture over $\rationals T_f$. If $f$ is
 injective, then $T_f$ fulfills the strong Atiyah conjecture even over
 $K T_f$ and we don't need 
the assumption $G\in\RAgroups$.
\end{example}
\begin{proof}
  Let $E_f$ be the direct limit of $G\stackrel{f}{\to }
G\stackrel{f}{\to} G\to\dots$. Then $E_f$ fulfills the strong
Atiyah conjecture by Proposition \ref{limits} or, if $f$ is
injective, by Proposition \ref{unions}. The claim follows from Proposition
\ref{amext} because $T_f$ is an extension of $E_f$
with quotient $\integers$.
\end{proof}

\section{Proofs of the Atiyah conjecture}\label{sec:proofs}

In this section we prove Theorem \ref{strongAt} and Proposition
\ref{classprop} and give additional information about
$\extendedC$.

\begin{lemma}\label{finextF}
  Let $H=F_1\times F_2\times \cdots$ be a direct sum of
  non-abelian free
  groups. Let $H\normalsubgroup G$ be a torsion-free finite extension
  of $H$. Then every finitely generated subgroup  of $G$ is contained
  in a group $V$ which fits into an
  extension $1\to U\to V\to A\to 1$ where $U$ is a subgroup of a
  direct sum of free groups and $A$ is torsion-free elementary amenable.
\end{lemma}
\begin{proof}
  By \cite[13.2]{Linnell(1998)} every finitely generated subgroup of
  $G$ is contained in a finite extension of a finitely generated
  product of free groups. So we may assume that $G$ (and 
  $H$) are finitely generated. We now use the ideas of
  \cite[13.4]{Linnell(1998)} where Linnell proves a related
  result. For a group $V$ let $S_nV$ be the
  intersection of normal subgroups of index $\le n$. Note that
  this is a characteristic subgroup, and if $V$ is finitely generated,
  $S_nV$ has finite index in $S$ by \cite[8.2.7]{Rowen(1988b)}. Set
  $H_n:= S_nF_1\times S_n F_2\dots S_n F_k$. Let $H_n'$ be the
  commutator subgroup of $H_n$. Then $H/H_n'$ is virtually
  abelian and
  torsion-free: it is the direct sum of $F_i/(S_nF_i)'$ containing the
  abelian subgroup $S_nF_i/(S_nF_i)'$ of finite index, and
  whenever $F$ is free and $N\normalsubgroup F$ is a normal subgroup
  then $F/N'$ is torsion free (compare e.g.~the last paragraph in the
  proof of 3.4.8 of \cite{Wolf(1984)}).

  Since $G$ is torsion-free, by \cite[13.4]{Linnell(1998)} for every
  prime number $p$ there is a
  subgroup $K_p\subgroup H$ which is normal in $G$ such that $G/K_p$ is
  virtually abelian and does not contain a finite subgroup of index
  $p$. For $n_p$ sufficiently
  large $H_{n_p}'\subgroup K_p$. We do this for all prime factors of
  $\abs{G/H}$ and let $n$ be the
  maximum of the $n_p$ (which exists since $G/H$ is finite). Then
  $1\to H_n'\to G\to A\to 1$ is an extension where $A$ is virtually
  abelian. We also have an extension $1\to H/H_n'\to A\to G/H\to
  1$. Since $H/H_n'$ is torsion-free, every torsion element of $A$ is
  mapped to an element of the same order in $G/H$. Therefore it
  suffices to check for a prime factor $p$ of $\abs{G/H}$ and $x\in
  A$ with $x^p=1$ that $x=1$. But we also have the extension $1\to
  K_p/H_n'\to A\to G/K_p\to 1$ where $K_p/H_n'\subgroup
  H/H_n'$ is torsion-free and $G/K_p$ does not contain
  $p$-torsion. This concludes the proof.
\end{proof}

\begin{proof}[Proof of Theorem \ref{strongAt} and Corollary
  \ref{corol:zero_divisors_in_D}]\strut\\
  Remember that every group in $\extendedC$ is torsion-free and
  belongs 
  by Lemma \ref{LinnelsC_in_RA} to $\RAgroups$.
  The result follows immediately from Proposition \ref{subgroups},
  Proposition \ref{limits}, and Proposition \ref{amext}.

  The corollary is an immediate consequence of the theorem (because of
  Lemma \ref{skew}). 
\end{proof}

\begin{proof}[Proof of Proposition \ref{classprop}]\strut\\
  By \cite{Magnus(1935)}, every finitely generated free group is residually torsion-free 
  nilpotent. Consequently, the same is true for every direct sum of
  free groups. Therefore, these groups are subgroups of inverse limits 
  of torsion-free nilpotent (and therefore elementary amenable)
  groups. It follows that direct sums of arbitrary free groups (as
 directed  unions of sums of finitely generated ones) and
  subgroups of these belong to $\extendedC$. By Lemma \ref{finextF},
  every torsion-free finite extension of a direct sum of free groups
  then also belongs to $\extendedC$. It now follows immediately from
  the definition of $\LinnelsC'$ that every torsion-free element of
  $\LinnelsC'$ belongs to $\extendedC$.

 Consider (for every ordinal $\alpha$) the subclasses
  $\LinnelsC_\alpha$ of $\LinnelsC$ with
  \begin{equation*}
\LinnelsC_0:=\freegroups\finiteGroups\quad\text{and}\quad
  \LinnelsC_{\alpha+1}:=(L\LinnelsC)\virtabelianGroups.
\end{equation*}
If $\alpha$ is
  a limit ordinal, set
  $\LinnelsC_\alpha:=\bigcup_{\beta<\alpha}\LinnelsC_\beta$. By
  induction, we proof that every
  torsion-free element of $\LinnelsC_\alpha$ belongs to
  $\extendedC$. We have just shown this for $\alpha=0$. If $\alpha$ is 
  a limit ordinal, nothing is to proof. By \cite[Lemma
  4.9]{Linnel(1993)}, each $L\LinnelsC_\alpha$ 
  is closed under extensions with finite quotients. If
  $G\in\LinnelsC_{\alpha+1}$ is torsion-free, we have an extension
  $1\to H\to G\xrightarrow{p}
  A\to 1$ with $A$ elementary amenable (since $A$ is abelian by
  finite) and $H\in L\LinnelsC_\alpha$. If $E\subgroup A$ is finite,
  therefore $p^{-1}(E)\in L\LinnelsC_\alpha$ (as a finite extension of 
  $H$). Moreover, $p^{-1}(E)$ is torsion-free, since it is a subgroup
  of the torsion-free group $G$. By induction,
  $p^{-1}(E)\in\extendedC$. It follows that $G\in\extendedC$. This
  finishes the induction step. Since
  $\LinnelsC=\bigcup_{\alpha}\LinnelsC_\alpha$ (by \cite[Lemma
  4.9]{Linnel(1993)}), the statement about $\LinnelsC$ follows.

  Next, we prove that $\extendedC$ is closed under direct sums. Since
  every direct sum is a union of finite direct sums, by induction it
  suffices to consider direct sums with two summands. Fix $G\in\extendedC$.

  The start of the induction is the trivial observation that
  $\{1\}\times G\in\extendedC$.
  For the induction step, we have to check that $H\times
  G\in\extendedC$ in the following three cases:
  \begin{enumerate}
  \item There is an exact sequence $1\to U\to H\xrightarrow{p} A\to 1$
    with $A$ 
    elementary amenable and such that $p^{-1}(E)\times G\in\extendedC$ 
    for every finite subgroup $E$ of $A$. But then we have an exact
    sequence 
    \begin{equation*}
      1\to U\times G\to H\times G\xrightarrow{q}A\to 1
    \end{equation*}
    and $q^{-1}(E)=p^{-1}(E)\times G$ for every finite subgroup $E$ of 
    $A$, and it follows that $H\times G\in\extendedC$.
  \item $H$ is a direct or inverse limit of a directed system of
    groups $H_i$, and $H_i\times G\in\extendedC$ for every $i$. Then
    $H\times G$ is the direct or inverse limit of the directed system
    $H_i\times G$ by Lemma \ref{lim1}, and therefore $H\times
    G\in\extendedC$.
  \item If $H\subset U$ and $U\times G\in\extendedC$, then $H\times
    G\subset U\times G$ which implies $H\times G\in\extendedC$.
  \end{enumerate}
  The detailed construction of families $\extendedC_\alpha$ such that
  this proof applies is standard (compare the similar construction of
  $\LinnelsC_\alpha$ used above) and is left to the reader.

  Every direct product is the inverse limit of finite direct sums. It
  follows that $\extendedC$ is also closed under taking direct
  products.

  We now proof by induction that $\extendedC$ is closed under free
  products with (arbitrary) free groups. The induction starts with the 
  observation that $\{1\}*F\in\extendedC$ if $F$ is free, what we have 
  already checked.

  For the induction step, we again have to check three cases.
  \begin{enumerate}
  \item Assume $H\subset U$ and $U*F\in\extendedC$ for every free
    group $F$. Then $H*F\subset U*F$ and consequently
    $H*F\in\extendedC$.
  \item Assume we have an exact sequence $1\to U\to H\to
    A\xrightarrow{p}1$ with $A$ elementary amenable and
    $p^{-1}(E)*F\in\extendedC$ for every free group $F$ and every
    finite subgroup $E$ of $A$. Consider the projection $q:=p*1\colon
    H*F\to A$. Then for every $E\subgroup A$ Lemma
    \ref{lem:free_product_inverse_image} implies $q^{-1}(E)\iso
    p^{-1}(E)*F'$ with some possibly
    different free group $F'$. Since the induction hypothesis applies
    to $p^{-1}(E)*F'$ if $E$ is finite, it follows that
    $H*F\in\extendedC$.
  \item Assume $H$ is the direct or inverse limit of a directed system 
    of groups $H_i$, and $H_i*F\in\extendedC$ for every free group
    $F$. By \cite[Lemma 2.9 and 2.10]{Schick(1998a)} $H*F$ is
    contained in the direct or inverse limit of the directed system
    $H_i*F$. It follows that $H*F\in\extendedC$.
  \end{enumerate}
  For arbitrary free products observe that
 $G*H$ is a subgroup of $(G\times H)*(G\times H)$ which is a subgroup
  of $(G\times H)*\integers$ (it is contained in the kernel of the
  projection onto $\integers$ by \cite[Lemma 5.5]{Gruenberg(1957)}).
  Since $\extendedC$ is subgroup closed, $G*H\in\extendedC$ if
  $G,H\in\extendedC$ follows from what we have proved about direct
  sums and free products with free groups. The class $\extendedC$ is
  closed under arbitrary free products by induction and considering
  directed unions (an arbitrary free product is the directed union of
  finite free products). This finishes the proof.
\end{proof}

The following lemmas were  needed in the proof of Proposition
\ref{classprop}.
\begin{lemma}\label{lim1}
  If $\pi$ is the direct or inverse limit of a system of groups
  $\pi_i$ and if $G$
  is any group, then
  $\pi\times G$ is the direct or inverse limit of $\pi_i\times G$.
\end{lemma}
\begin{proof}
  Compare e.g.~in \cite[2.8 and 2.9]{Schick(1998a)}.
\end{proof}

\begin{lemma}\label{lem:free_product_inverse_image}
  Assume $p\colon H\to A$ a homomorphism of groups. Consider
  $q:=p*1\colon H*G\to A$. Then, for every subgroup $E$ of $A$,
  $q^{-1}(E)= p^{-1}(E)* (*_{t\in T} (t^{-1}Gt))$, where $T$ is a
  transversal for $p^{-1}(E)$ in $H$, i.e.~a system of representatives 
  for the cosets $H/p^{-1}(E)$.
\end{lemma}
\begin{proof}
  Obviously, $q$ maps $p^{-1}(E)$ and $t^{-1}Gt$ to $E$. On the other
  hand, a normal form argument implies that each element in
  $q^{-1}(E)$ lies in the subgroup of $H*G$ generated by $p^{-1}(E)$ and
  $t^{-1}Gt$ ($t\in T$). The normal form argument also implies that
  this subgroup is isomorphic to the free product $p^{-1}(E)* (*_{t\in
    T}(t^{-1}Gt))$. 
\end{proof}



We now explain how to fill the gap (pointed out in Section
\ref{sec:introduction}) of Linnell's
result about the Atiyah conjecture for the class $\LinnelsC'$, using
our methods (and obtaining a slightly weaker conclusion).

\begin{corollary}\label{corol:LinnelsCprime}
  Assume $G\in\LinnelsC'$ and there
  is a bound on the orders of the finite subgroups of $G$. Then the
  strong Atiyah conjecture is true over $\rationals G$.
\end{corollary}
\begin{proof}
  Linnell's proof in \cite[Section 13]{Linnell(1998)} applies as soon
  as the strong Atiyah conjecture over $\rationals G$ is established
  for every  direct sum
  of free groups $G$. However, we have seen that these groups belong
  to $\extendedC$ and therefore satisfy the strong Atiyah conjecture
  over $\rationals G$.
\end{proof}

\section{Generalizations}\label{sec:gen}

In the second part of the paper we concentrated on torsion-free
groups. By Linnell's methods, many of the results hold if instead of
this one requires a bound on the orders of finite subgroups. In particular,
a version of Proposition \ref{amext}  remain true if
$Q$ is not
  assumed to be torsion-free but only has a bound on the orders of its
  finite subgroups. One proves this by transfinite induction, using a
  little bit more of the machinery of Linnell
  \cite{Linnel(1993),Linnell(1998)}. Correspondingly, one can then
  (similarly to $\extendedC$)
   define a larger class of groups (not necessarily torsion-free)
  for which the strong Atiyah
  conjecture is true.

Another possible generalization is given by only considering the
$1$-dimen\-sio\-nal Atiyah conjecture. This is defined as in Definition
\ref{AtiyahCon}, but one does consider only $1\times 1$-matrices $A$,
i.e.~$A\in K G$. This is interesting because already the strong
$1$-dimensional Atiyah conjecture for
torsion-free groups implies the zero divisor conjecture. Moreover,
Linnell \cite{Linnel(1992)} proves that the $1$-dimensional strong
Atiyah conjecture for
torsion-free groups is stable under extensions with right-orderable
groups.

One immediately checks that there is a $1$-dimensional versions of
Proposition 
\ref{limits}. Hence we can define a class of groups similar to
$\extendedC$ but which is
closed under extensions with
right-orderable quotient (instead of extensions with elementary
amenable quotient), and the $1$-dimensional strong Atiyah conjecture
over $\rationals G$ is true if $G$ is
contained in this class. In
particular, $\rationals G$ is zero divisor free. 

\section{Treelike groups}\label{treechar}

\begin{proposition}\label{gentrees}
  Let $H$ be a finite group. It follows that $H$ fulfills  the Atiyah
  conjecture of order
  $\frac{1}{L}\integers$ over $KH$, with $L=\abs{H}$. Let $G$ be a
  group and $\Omega$ and $\Delta$
  sets with commuting $G$-action from the left and $H$-action from the
  right such that $\Omega$ and $\Delta$ are free $H$-sets. Let
  $\Omega=\Omega'\disjointunion X$ and assume that $\Delta$ and $\Omega'$ are free
  $G$-sets and $X$ consists of $1\le r<\infty$ $H$-orbits. Assume that
  there is a bijective map $\phi:\Delta\to\Omega$ which is an
  $H$-almost $G$-map and a right $H$-blockwise map (compare Definition
  \ref{defnota}).

Then $G$
  fulfills the Atiyah conjecture of order
  $\frac{1}{rL}\integers$ over $KG$.
\end{proposition}
We used the following notation:
\begin{definition}\label{defnota}
  The map $\phi:\Delta\to\Omega$ is an \emph{$H$-almost $G$-map} if for each
  $g\in G$ the set $\{x\in\Delta\|\; \phi(gx)\ne
  g\phi(x)\}$ is contained in only finitely many $H$-orbits.\\
  It is called a \emph{right $H$-blockwise} map if each $H$-orbit of
  $\Delta$ is mapped bijectively to an
  $H$-orbit of $\Omega$.
\end{definition}

\begin{remark}
  Whenever the assumptions of Proposition \ref{gentrees} are
  fulfilled, we can replace $H$ by $\{1\}$ (and therefore get
  $L=1$). However, $X$ then consists of $r\times\abs{H}=rL$
  $\{1\}$-orbits (i.e.~points), and the conclusion is unchanged.
\end{remark}

\begin{proof}[Proof of Proposition \ref{gentrees}]
 Our proof generalizes Linnell's approach in
  \cite{Linnell(1998)}.

  Assume that we have the $G$-$H$ sets $\Omega$ and $\Delta$ as above.
  Given $A\in M_n(K G)$, the $G$-operation induces bounded operators
  $A_\Delta$ on $l^2(\Delta)^n$ and $A_{\Omega'}$ on $l^2(\Omega')^n$ in
  the obvious way. Set $A_\Omega:=A_{\Omega'}\oplus 0$ on
  $l^2(\Omega)^n=l^2(\Omega')^n\oplus l^2(X)^n$. 

  Let $P\in M_n(\NeumannN G)$ be the projection onto the image of $A$.
  Since $\Delta$ and
  $\Omega'$ are free $G$-sets one can extend the above construction to
   $P$. It follows immediately that $P_\Delta$ is the projection onto
  the image of $A_\Delta$ and $P_\Omega$ the projection onto the image
  of $A_\Omega$. Since
  the $H$-actions commute with the $G$-actions, all operators
  constructed in this way are $H$-equivariant.

 There are only finitely many $g\in G$ such that $A=:(a_{ij})\in
M_n(\complexs G)$ has a non-trivial coefficient of $g$ in at least one of
the $a_{ij}$. For each of these, there are only finitely many
$H$-orbits in $\Delta$ on which $\phi(gx)\ne g\phi(x)$. Let
$\Delta_0$ be the union of
  this finite number of $H$-orbits in $\Delta$, and set
$\Delta_c=\Delta-\Delta_0$ and $\Omega_c:=\phi(\Delta_c)$.

Restricted to $l^2(\Delta_c)^n$, $\phi^* A_\Omega\phi$ and $ A_\Delta$
coincide:
\begin{equation}
  \label{eq:coincide}
  A_\Delta |_{l^2(\Delta_c)^n} = \phi^*A_\Omega\phi |_{l^2(\Delta_c)^n}.
\end{equation}

Choose an exhaustion $\Delta_0\subset \Delta_1\subset
  \Delta_2\subset\dots\Delta$ of $\Delta$
  (i.e.~$\Delta=\bigcup_{k\in\naturals} \Delta_k$) where each $\Delta_k$
  consists of finitely many (entire) $H$-orbits. Then (since
  $\phi$ is an $H$-blockwise map)
  $\Omega_k:=\phi(\Delta_k)$ also consists of finitely many 
  $H$-orbits.

  If $\Delta$ consists of  uncountably many $H$-orbits,  we have to
  use an exhaustion index by an uncountable directed index set. The
  proof remains essentially unchanged. Since our applications only
  deal with the countable setting for simplicity we
  stick to the index set $\naturals$.

As $\Delta=\bigcup_{k\in\naturals}\Delta_k$ we have
$A_\Delta(l^2(\Delta)^n)= \overline{\bigcup_k
  A_\Delta(l^2(\Delta_k)^n)}$.

 Let $P_{\Delta,k}$ be the projection onto
$\overline{A_\Delta(l^2(\Delta_k)^n)}$, and
$P_{\Omega,k}$ the projection onto
$\overline{A_\Omega(l^2(\Omega_k)^n)}$. 

Because in the matrix $A=(a_{ij})$ only finitely many elements of $G$ have
non-trivial coefficient in at least one of the $a_{ij}$, $A$
has finite
propagation, i.e.~for every $k\in \naturals$ we find $n(k)\in\naturals$
such that $A_\Delta(l^2(\Delta_k)^n)\subset
l^2(\Delta_{n(k)})^n$. Similarly (if $n(k)$ is chosen big enough)
$A_\Omega(l^2(\Omega_k)^n)\subset l^2(\Omega_{n(k)})^n$. Particularly, $A_\Delta(l^2(\Delta_k)^n)$ and
$A_\Omega(l^2(\Omega_k)^n)$ are the ranges of finite matrices over
$K H$ to which the Atiyah conjecture applies and by assumption on $H$ 
\begin{equation}\label{eq:k-rk}
  \begin{split}
    \tr_{H}(P_{\Delta,k}) &= \dim_{ H}(A_\Delta(l^2
    \Delta_k)^n)\in \frac{1}{L}\integers,\\
    \tr_{H}(P_{\Omega,k}) &= \dim_{H}(A_\Omega(l^2
    \Omega_k)^n)\in \frac{1}{L}\integers.
\end{split}
\end{equation}

Let $P_{\Delta,k,c}$ be the projection onto
$A_\Delta(l^2(\Delta_k\cap\Delta_c)^n)$ and $P_{\Delta,c}$ the
projection onto $A_\Delta(l^2(\Delta_c)^n)$. Define $P_{\Omega,k,c}$ and
$P_{\Omega,c}$ correspondingly. By Equation \eqref{eq:coincide}
$P_{\Delta,k,c}= \phi^*P_{\Omega,k,c}\phi$ and
$P_{\Delta,c}=\phi^*P_{\Omega,c}\phi$. In particular, if
$x\in\Delta^n\subset l^2(\Delta)^n$ we have (with $\phi$ diagonally
extended to $\Delta^n$)
\begin{equation}\label{eq:vanish}
  \begin{split}
    \innerprod{P_{\Delta,k,c}x,x}_{l^2(\Delta)^n} &=\innerprod {
      P_{\Omega,k,c}\phi(x),\phi(x)}_{l^2(\Omega)^n},\\
    \innerprod{P_{\Delta,c}x,x}_{l^2(\Delta)^n} &=\innerprod {
      P_{\Omega,c}\phi(x),\phi(x)}_{l^2(\Omega)^n}.
\end{split}
\end{equation}
Set $Q_{\Delta,k}:= P_{\Delta,k}-P_{\Delta,k,c}$,
$Q_\Delta:=P_\Delta-P_{\Delta,c}$, and define $Q_{\Omega,k}$ and
$Q_\Omega$ correspondingly. These are the projections onto the
complement of the smaller sets in the larger ones. 

Note that the
images of the $Q$'s are already obtained from $l^2(\Delta_0)^n$ and
$l^2(\Omega_0)^n$,
respectively. More precisely, let $A_0$ be the restriction of $A_\Delta$ to
$l^2(\Delta_0)^n$. Then $Q_{\Delta,k} A_0$ induces a weakly exact
sequence (i.e.~the images are dense in the kernels)
\begin{equation}\label{eq:sur1}
0\to\ker(Q_{\Delta,k} A_0)\to l^2(\Delta_0)^n\xrightarrow{Q_{\Delta,k}
A_0}\im(Q_{\Delta,k})\to 0.
\end{equation}
We only have to check weak exactness at $\im(Q_{\Delta,k})$. 
 $Q_{\Delta,k}(A(l^2(\Delta)^k))$ is dense in $\im(Q_{\Delta,k})$. But
 $A(l^2(\Delta_k)^n)=
 A(l^2(\Delta_0)^n)+A(l^2(\Delta_k\cap\Delta_c)^n)$ and
 $Q_{\Delta,k}(A(l^2(\Delta_k\cap\Delta_c)^n))=0$. 
 For the kernel
we have  the weakly exact
sequence 
\begin{multline}\label{eq:split1}
0\to\ker(A_0)\to \ker(Q_{\Delta,k}
A_0) \stackrel{A_0}{\to}\\
\overline{A_\Delta(l^2(\Delta_0)^n)\cap
A_\Delta(l^2(\Delta_k\cap\Delta_c)^n)}\to 0,
\end{multline}
which would be clear with $A_\Delta(l^2(\Delta_k\cap\Delta_c)^n)$ replaced
by $\ker(Q_{\Delta,k})$ ---but
 if $f\in A_\Delta(l^2(\Delta_k)^n)$ and in particular if
 $f\in A_\Delta(l^2(\Delta_0)^n)$, then $f\in\ker(Q_{\Delta,k})$ if and
 only if
$f\in A_\Delta(l^2(\Delta_k\cap\Delta_c)^n)$. This is true since all the
vector spaces we are considering here are finite dimensional (a
consequence of the assumption that $H$ is finite).

In the same way one gets the weakly exact sequences
\begin{gather}\label{eq:sur2}
0\to \ker(Q_\Delta A_0)\to
l^2(\Delta_0)\xrightarrow{Q_\Delta A_0} \im(Q_\Delta)\to 0,\\
  \label{eq:split2}
   0 \to \ker(A_0) \to \ker(Q_\Delta A_0) \stackrel{A_0}{\to}
  \overline{A_\Delta(l^2(\Delta_0)^n)\cap
A_\Delta(l^2(\Delta_c)^n)} \to 0.
\end{gather}
Observe that 
\begin{equation*}
  \overline{A_\Delta(l^2(\Delta_0)^n) \cap
A_\Delta(l^2(\Delta_c)^n)} = \overline{\bigcup_{k\in\naturals} {A_\Delta(l^2(\Delta_0)^n)}\cap
{A_\Delta(l^2(\Delta_k\cap \Delta_c)^n)}},
\end{equation*}
hence because of continuity of $\dim_H$ \cite[1.4]{Lueck(1997)}
\begin{multline}\label{conv}
  \dim_H(\overline{A_\Delta(l^2(\Delta_0)^n) \cap A_\Delta(l^2(\Delta_c)^n)}) \\
  =
    \lim_{k\to\infty} \dim_H(\overline{A_\Delta(l^2(\Delta_0)^n)\cap
    A_\Delta(l^2(\Delta_k\cap \Delta_c)^n)}).
\end{multline}
Using the additivity of $\dim_H$ under short weakly exact sequences \cite[1.4]{Lueck(1997)},
Equation \eqref{conv}
implies together with Equation \eqref{eq:sur1}, \eqref{eq:split1}, \eqref{eq:sur2},
 and \eqref{eq:split2}
(where all numbers are finite and bounded by $\dim_H(l^2(\Delta_0)^n)$)
\begin{multline}
  \label{eq:conv}
  \tr_H(Q_\Delta)= \dim_H(\im(Q_\Delta)) = \lim_{k\to\infty}
  \dim_H(\im(Q_{\Delta,k}))\\
   =\lim_{k\to\infty}
  \tr_H(Q_{\Delta,k}). 
\end{multline}

Exactly in the same way one obtains
\begin{equation}
  \label{eq:convO}
  \tr_H(Q_\Omega)= \lim_{k\to\infty}
  \tr_H(Q_{\Omega,k}). 
\end{equation}

Choose a set of representatives $T$ for the $H$-orbits of
$\Delta$. Because $\phi$ is an $H$-blockwise bijective map, $\phi(T)$
will be a set of representatives for the $H$-orbits of $\Omega$. For
$t\in T$ set
$t_i:=(0,\dots,0,t,0,\dots,0)\in \Delta^n$ ($t$ at the $i$-th
position). Then
(by the definition of the $H$-trace)
\begin{equation*}
  \begin{split}
    \tr_H(Q_\Delta)= & \sum_{t\in T}\sum_{i=1}^n \innerprod{Q_\Delta
      t_i,t_i}_{l^2(\Delta)}, \\
    \tr_H(Q_\Omega) = &  \sum_{t\in T}\sum_{i=1}^n\innerprod{Q_\Omega
      \phi(t_i),\phi(t_i)}_{l^2(\Delta)}.
\end{split}
\end{equation*}

Because $\Delta$ is a free $G$-set, $G\cdot
t\subset\Delta$ can for each $t\in T$ be identified with $G$ and therefore 
\begin{equation*}
  \sum_{i=1}^n \innerprod{P_\Delta t_i,t_i}_{l^2(\Delta)^n} =
  \tr_G(P)=\dim_G(\im(A)).
\end{equation*}
Similarly, since $\Omega'$ is a free $G$-set, if $\phi(t)\in\Omega'$
then
\begin{multline*}
    \sum_{i=1}^n\innerprod{P_\Omega \phi(t_i),\phi(t_i)}_{l^2(\Omega)^n} =
    \sum_{i=1}^n\innerprod{P_\Omega
    \phi(t_i),\phi(t_i)}_{l^2(\Omega')^n}=\\
  \tr_G(P)=\dim_G(\im(A)).
  \end{multline*}
However, if $\phi(t)\in X$ then by the construction of $P_\Omega$
\begin{equation*}
    \innerprod{P_\Omega \phi(t_i),\phi(t_i)}_{l^2(\Omega)} = 0.
\end{equation*}

Therefore 
\begin{equation}\label{classic}
  \sum_{t\in T}\sum_{i=1}^n (\innerprod{P_\Delta
  t_i,t_i}_{l^2(\Delta)^n} - \innerprod{P_\Omega
  \phi(t_i),\phi(t_i)}_{l^2(\Omega)^n}) = r \dim_G(\im(A))
\end{equation}
(where $r=\abs{\phi(T)\cap X}$ is the number of $H$-orbits in
$X$). Note that the sum is finite so that convergence is not an
issue here. We want to use this equation to determine the value of
$\dim_G(\im(A))$. 

From our splitting of $P_\Delta$ and $P_\Omega$ we see that for each
$t\in T$ and  $i\in\{1,\dots,n\}$
\begin{equation}\label{inf_eq}
  \begin{split}
    &\innerprod{P_\Delta t_i,t_i}_{l^2(\Delta)^n} - \innerprod{P_\Omega
      \phi(t_i),\phi(t_i)}_{l^2(\Omega)^n}\\
    &= \underbrace{\innerprod{P_{\Delta,c}t_i,t_i} -
     \innerprod{P_{\Omega,c}\phi(t_i),
      \phi(t_i)}}_{\stackrel{\eqref{eq:vanish}}{=} 0} +
      \innerprod{Q_\Delta t_i,t_i} - \innerprod{Q_\Omega\phi(t_i),\phi(t_i)}.
\end{split}
\end{equation}
In the same way
\begin{equation}\label{k_eq}
   \begin{split}
    &\innerprod{P_{\Delta,k} t_i,t_i}_{l^2(\Delta)^n} - \innerprod{P_{\Omega,k}
      \phi(t_i),\phi(t_i)}_{l^2(\Omega)^n}\\
    &= \underbrace{\innerprod{P_{\Delta,k,c}t_i,t_i} -
      \innerprod{P_{\Omega,k,c}\phi(t_i),
      \phi(t_i)}}_{\stackrel{\eqref{eq:vanish}}{=} 0} +
      \innerprod{Q_{\Delta,k} t_i,t_i} -
      \innerprod{Q_{\Omega,k}\phi(t_i),\phi(t_i)}.
\end{split}
\end{equation}
Note that by Equation \eqref{eq:k-rk} for each finite $k$, all the operators in
Equation \eqref{k_eq} are
of finite $H$-rank individually and  therefore we can write
\begin{equation}
  \label{eq:rkok}
  \begin{split}
    \tr_H(Q_{\Delta,k})-\tr_H(Q_{\Omega,k}) & = \sum_{t\in
      T}\sum_{i=1}^n \innerprod{Q_{\Delta,k} t_i,t_i} -
    \innerprod{Q_{\Omega,k}\phi(t_i),\phi(t_i)} \\
   & \stackrel{\eqref{k_eq}}{=} \tr_H(P_{\Delta,k}) -
    \tr_H(P_{\Omega,k}) \stackrel{\eqref{eq:k-rk}}{\in}
    \frac{1}{L}\integers .
\end{split}
\end{equation}

This implies, since $\frac{1}{L}\integers$ is discrete,
\begin{equation}
  \label{eq:10}
    \tr_H(Q_\Delta)-\tr_H(Q_\Omega)
    \stackrel{\text{\eqref{eq:conv}, \eqref{eq:convO}}}{=}
    \lim_{k\to\infty}
    \left( \tr_H(Q_{\Delta,k})- \tr_H(Q_{\Omega,k})\right)
    \stackrel{\eqref{eq:rkok}}{\in} \frac{1}{L}\integers .
  \end{equation}
Putting everything together we see
\begin{equation*}
  \begin{split}
    r\dim_H(\im(A)) & \stackrel{\eqref{classic}}{=} \sum_{t\in T}\sum_{i=1}^n (\innerprod{P_\Delta
      t_i,t_i}_{l^2(\Delta)^n} - \innerprod{P_\Omega
      \phi(t_i),\phi(t_i)}_{l^2(\Omega)^n})\\
     & \stackrel{\eqref{inf_eq}}{=} \sum_{t\in T}\sum_{i=1}^n (\innerprod{Q_{\Delta} t_i,t_i} -
     \innerprod{Q_{\Omega}\phi(t_i),\phi(t_i)})\\
     & =
    \tr_H(Q_\Delta)-\tr_H(Q_\Omega) \stackrel{\eqref{eq:10}}{\in} \frac{1}{L}\integers.
  \end{split}
\end{equation*}
Since $A\in M_n(\complexs G)$ was arbitrary, Proposition \ref{gentrees}
follows.
\end{proof}

\begin{definition}
  A group $G$ is called {\em treelike} (of {\em finity} $r$), if there are
  free $G$-sets $\Delta$ and $\Omega'$, a non-empty trivial $G$-set $X$ with 
$0\ne r\in \naturals$ elements and an almost $G$-map 
\begin{equation*}
\alpha: \Delta \to
\Omega' \disjointunion X .
\end{equation*}
\end{definition}

Julg and Valette \cite{Julg-Valette(1984)}  show that free groups are
 treelike of finity $1$. On the other hand, Dicks 
and Kropholler \cite{Dicks-Kropholler(1995)}  prove that if $G$ is
 treelike of finity $1$ then $G$ is free.

As an immediate consequence of Proposition \ref{gentrees} with
  $H=\{1\}$ we get:
\begin{theorem}\label{treeext}
  Let $Q$
be a  treelike group of finity $r$. Then 
$Q$  fulfills the  Atiyah conjecture of order
  $\frac{1}{r}\integers$ over $\complexs Q$. 
\end{theorem}

We generalize in this section the main result in
\cite{Dicks-Kropholler(1995)} and show  that the
class of  treelike groups of finity $r$ coincides with
the class of groups which contain a  free group of index $r$. So, in view of
Lemma \ref{subgroups},
Theorem \ref{treeext} does not provide new cases of the Atiyah
conjecture.

\begin{theorem}\label{characttree}
  The following statements are equivalent for a group $G$:
  \begin{enumerate}
  \item\label{treelike} $G$ is treelike of finity which divides $r$.
  \item\label{treeact} $G$ acts on a tree and all vertex stabilizers
    have order which
    divides $r$.
  \item\label{graphgroups} $G$ is the fundamental group of a graph of
    groups where the
    order of each vertex group divides $r$.
  \item\label{virtfree} $G$ contains a free subgroup of finite index
    $d$ which
    divides $r$. 
  \item\label{forest} $G$ acts freely on a union of $d$ trees where
    $d$ divides $r$.
  \item\label{freefin} $G$ contains a free subgroup of finite index and the order of
    every finite subgroup of $G$ divides $r$.
  \end{enumerate}
\end{theorem}
 All
  statements which don't involve treelikeness are well
  known. Those proofs for which no reference in a
  standard textbook could be found are included here.
\begin{proof}
  \ref{treelike} $\implies$ \ref{treeact}: 
 Suppose $G$ is treelike of finity $r$. That means we
  have two free $G$-sets $\Delta$ and $\Omega'$ and
  a bijective almost $G$-map $\theta:\Delta\to\Omega:=\Omega'\amalg
  (\amalg_{i=1}^r\{*\})$. Let $V$ be the set of all maps with same
  domain and range which coincide with $\theta$ outside a finite
  subset of $\Delta$. By \cite[2.4]{Dicks-Kropholler(1995)} $V$ is the
  vertex set of a $G$-tree with finite edge stabilizers. We want to
  show that the vertex stabilizers are finite, and their order is
 divisible by $r$. 

  That means we want to show that there
  is no way to change $\theta$ at only finitely many points in the
  domain to obtain an $f$ fixed by $H\subgroup G$ unless $H$ is finite and its
  order divides $r$. Remember that the action of $G$ on
  the maps from $\Delta$ to $\Omega$ is given by conjugation:
  $f^g(x)=g^{-1}f(gx)$. Therefore, $f$ being an $H$-fixed point  is equivalent
  to $f$ being an $H$-map. 

Each  map $f\in V$ has an index
  $i(f):= \sum_{x\in\Omega}( \abs{f^{-1}(x)}-1)$. This is well defined
  if the inverse image of almost every point of $\Omega$ consists of
  one point. In particular, it is well defined for the bijective map
  $\theta$ and $i(\theta)=0$. If we change a map at one point (or
  iteratively at finitely many points), this
  index remains well defined and unchanged (the inverse image of one
  point will be smaller by one, whereas the inverse image of another
  point will be larger by one). In particular, $i(f)=0$ for every
  $f\in V$.

 Assume that $f\in V$ is an $H$-map for
  $H\subgroup G$. Then $f$ will miss, say, $k$ of the $r$ $H$-fixed points of
  $\Omega$, and the inverse image of the remaining $r-k$ of those will
  consist of unions of free $H$-orbits.  In addition, the inverse image of each
  point in a given free $H$-orbit of $\Omega$ will be of equal
  size. If $H$ is infinite then $k=r$ because no point can have an infinite
  inverse image (else $i(f)$ is not defined). But then $f$ is not
  injective because $i(f)=0$. Therefore the inverse image of  at least
  one free $H$-orbit in $\Omega$ consists of more than one orbit,
  which is impossible because $i(f)$ is defined. This contradiction
  shows that no infinite subgroup $H$ of $G$ fixes a point $f\in V$.

  If $H$ is finite and $f$ is an $H$-map then the fact that $\Delta$
  and $\Omega'$ are free $H$-sets implies as above 
 $0=i(f) = -r   + \abs{H}\cdot N$, i.e.~$r$ is
  divisible by the order of $H$.

\ref{treeact} $\implies$ \ref{graphgroups}: This is the
  structure theorem for groups acting on a tree, compare
  e.g.~\cite[I.4.1.]{Dicks-Dunwoody(1989)}.

\ref{graphgroups} $\implies$ \ref{virtfree}
  This is proved as \cite[Theorem 8.3]{Bass(1993)}. The statement
   follows also easily from the arguments of
  \cite[IV.1.6]{Dicks-Dunwoody(1989)}, where is is shown that $G$
  contains a free normal subgroup of finite index. 


  \ref{virtfree}$\implies$\ref{forest}: Let $F$ be a free subgroup of
  index $d$ in $G$.
Choose a free basis $\{f_i\}$ of $F$. Construct the following
  graph: the elements of $G$ are the 
  vertices, and two vertices $g$ and $h$ are joined by an edge if and only if
  $gf_i=h$ for a suitable $f_i$ in the basis of $F$.

  The graph is a free left $G$-set in the obvious way (with transitive
  action on the vertices). The component of the trivial element is the
  Cayley graph of $F$ with respect to the generators $\{f_i\}$ and
  therefore is a tree \cite[I.8.2]{Dicks-Dunwoody(1989)}. The
  vertex sets of the other
  components are the translates of $F$. Consequently, the graph consists
  of $d$ trees.

\ref{forest} $\implies$ \ref{treelike}:
A construction of  Julg and Valette gives an almost $G$-map from
  the set of edges to the set of vertices united with $d$ additional
  points which are fixed under the $G$-action. It is done as follows:
  choose base points for each of the trees. These base points are mapped
  to the $d$ additional points. Each other vertex is mapped to the first
  edge of the geodesic starting at this vertex and ending at the
  basepoint of its
  component. This geodesic and therefore the edge is unique because we
  are dealing with trees. Since the $G$-action was free, we get a
  treelike structure for $G$ with finity $d$.

\ref{freefin} $\implies$ \ref{treeact}: A group which
 contains a free subgroup of finite index  acts on a tree with
finite vertex stabilizers \cite[IV.1.6]{Dicks-Dunwoody(1989)}. The
stabilizers are subgroups of $G$, therefore their order divides $r$. 

\ref{graphgroups} $\implies$ \ref{freefin}: By
\cite[IV.1.6]{Dicks-Dunwoody(1989)}, if $G$ is the fundamental group of
a graph of groups with finite vertex groups where the order is bounded
by $r$, then it contains a free subgroup of finite index. Moreover,
every finite subgroup of $G$ is contained in a conjugate of a vertex
group \cite[I.7.11]{Dicks-Dunwoody(1989)}, i.e.~the order of every
finite subgroup divides $r$.
\end{proof}



\end{document}


\pagestyle{myheadings}
\markboth{Thomas Schick}{Erratum concerning ``amenable quotients''}

\date{Last compiled \today}

\title{Errata to ``Integrality of $L^2$-Betti numbers'', ``$L^2$-determinant
  class and approximation of $L^2$-Betti nubmers'', and work based on these}

\author{Thomas Schick\thanks{
\protect\href{mailto:schick@uni-math.gwdg.de}{e-mail:
  schick@uni-math.gwdg.de}
\protect\\
\protect\href{http://www.uni-math.gwdg.de/schick}{www:~http://www.uni-math.gwdg.de/schick}
}\\
Georg-August-Universit\"at G{\"o}ttingen\\
Germany}
\maketitle

In \cite{Schick}*{Section 4}, a general notion of ``amenable
extension'' $U\subgroup \pi$ is defined, and in \cite{Schick}*{Definition
  1.12} a class of group $\mathcal{G}$ is defined which is in particular
closed under such (generalized) amenable extensions.

A particular example of a generalized amenable extension is a normal amenable
extension, i.e.~$U$ is a normal subgroup of $\pi$ and $\pi/U$ is amenable.

The main result, \cite{Schick}*{Theorem 1.14} is then proved for groups in
the class $\mathcal{G}$, based on claimed stability of the relevant properties
of groups under ``amenable extension''.

Unfortunately, as pointed out by Christian Wegner, the proof of the relevant
stability and approximation property for the general notion of amenable
extension is not valid,
it is based on a commutation relation which can't be established.

Nonetheless, the proof works perfectly well for usual normal amenable
extensions. Therefore, the assertions of the paper must be restricted to
normal amenable extensions; in particular, the definition of $\mathcal{G}$
must be modified such that ``amenable extension'' has to be replaced by the (a
priori more restrictive) notion of ``normal amenable extension''.

The definition of $\mathcal{G}$ and the approximation result for amenable
extensions has been taken up in \cite{MR1777117}*{Definition 3 and Proposition
  1}. Consequently, also in this paper each use of ``amenable extension'' has
to be restricted to ``normal amenable extension'', the class $\mathcal{G}$ has
to be redefined accordingly. 

The defintiion of $\mathcal{G}$ has also been taken up in
\cite{DLMSY}*{Definition 1.8}, and, based on the methods of
\cite{MR1777117} generalized approximation theorems and further properties are
claimed to be established for the groups in $\mathcal{G}$. Again, these
statements are established only if, throughout, ``amenable extension'' is
changed to ``normal amenable extension'' and the class $\mathcal{G}$ has to be
replaced by the (a priori smaller) class which is closed only under ``normal
amenable extension''. 

Similarly, the definition of $\mathcal{G}$ is taken up in
\cite{MR2200261}*{Definition 4.6}, and based on it a class $\hat{\mathcal{G}}$
is defined; and the results of \cite{MR1777117} are used. Therefore, as
before, the definitions have to be modified to allow only ``normal amenable
extensions'' to have valid proofs for the statements made about groups in
$\mathcal{G}$ in \cite{MR2200261}.


Similarly, the definition of $\mathcal{G}$ has been taken up in
\cite{Knebusch}*{Situation 3.1} and the notion of (generalized) ``amenable
extension'' in \cite{Knebusch}*{Definition 5}. The results stated in
\cite{Knebusch} for generalized ``amenable extensions'' and for groups in
$\mathcal{G}$, e.g.~ are generalizations of and based on the methods of
\cite{MR1777117}; consequently they again have to be modified by replacing the
(generalized) ``amenable extensions'' by ``normal amenable extensions''
throughout, and by using the (a priori) smaller class $\mathcal{G}$ based on
this.

\begin{remark}
  To the authors knowledge, no example of an (generalized) amenable extension
  which is not a normal amenable extension is know, in particular no such
  example  and has been used explicitly in the literature.
\end{remark}